\documentclass[12pt]{amsart}
\usepackage{times,epsf}

\begin{document}

\newtheorem{thm}{Theorem}[section]
\newtheorem{lem}[thm]{Lemma}
\newtheorem{cor}[thm]{Corollary}

\theoremstyle{definition}
\newtheorem{defn}{Definition}[section]

\theoremstyle{remark}
\newtheorem{rmk}{Remark}[section]

\def\square{\hfill${\vcenter{\vbox{\hrule height.4pt \hbox{\vrule
width.4pt height7pt \kern7pt \vrule width.4pt} \hrule height.4pt}}}$}
\def\T{\mathcal T}

\newenvironment{pf}{{\it Proof:}\quad}{\square \vskip 12pt}

\title{Mean Convex Hulls and Least Area Disks Spanning Extreme Curves}
\author{Baris Coskunuzer}
\address{Department of Mathematics \\ Yale University \\ New Haven CT 06520}
\email{baris.coskunuzer@yale.edu}

\maketitle

\newcommand{\cirD}{\overset{\circ}{D}}
\newcommand{\Si}{S^2_{\infty}({\Bbb H}^3)}
\newcommand{\PI}{\partial_{\infty}}
\newcommand{\BH}{\Bbb H}
\newcommand{\BR}{\Bbb R}
\newcommand{\BC}{\Bbb C}
\newcommand{\BZ}{\Bbb Z}

\begin{abstract}

We show that for any extreme curve in a $3$-manifold $M$, there exist a canonical mean convex hull containing
all least area disks spanning the curve. Similar result is true for asymptotic case in $\BH^3$ such that for
any asymptotic curve $\Gamma\subset \Si$, there is a canonical mean convex hull containing all minimal planes
spanning $\Gamma$. Applying this to quasi-Fuchsian manifolds, we show that for any quasi-Fuchsian manifold,
there exist a canonical mean convex core capturing all essential minimal surfaces. On the other hand, we also
show that for a generic $C^3$-smooth curve in the boundary of $C^3$-smooth mean convex domain in $\BR^3$,
there exist a unique least area disk spanning the curve.

\end{abstract}

\section{Introduction}

We study the Plateau problem for extreme curves in a Riemannian $3$-manifold. The existence of least area
disks for any simple closed curve in a Riemannian manifold was proved by Morrey half a century ago. The
regularity properties and number of solutions questions have been studied in the following decades. On the
other hand, for extreme curves, this problem has many interesting properties. Meeks-Yau, Hass-Scott and
Fanghua Lin have studied the Plateau problem in mean convex domains, and proved many important features of
these least area disks, [MY1], [HS], [Li].

In this paper, we will concentrate on the same problem. We improve and reformulate some known results in a
natural way, and try to clarify the picture with a simplified proof by using the techniques of [MY1] and
[HS]. Moreover, we will extend these properties to asymptotic Plateau problem in $\BH^3$. On the other hand,
we prove a new generic uniqueness result for extreme curves in $\BR^3$ by using some topological techniques
and analytical results of [TT].

Now, we list the main results of the paper. First, we show existence of a canonical object called {\em mean
convex hull} for any extreme curve.\\

\noindent \textbf{Theorem 3.2.} Let $\Omega$ be a mean convex domain in a $3$-manifold $M$ and
$\Gamma\subset\partial \Omega$ be a simple closed curve. Then either there exist a unique least area disk
$\Sigma$ in $\Omega$ with $\partial \Sigma = \Gamma$, or there exist a canonical mean convex hull $N$  in
$\Omega$ such that $\partial N= \Sigma^+ \cup \Sigma^-$ where $\Sigma^\pm$ are uniquely defined extremal
least area disks in $\Omega$ with $\partial \Sigma^\pm = \Gamma$. Moreover, all least area disks $\Sigma
'\subset \Omega$ spanning $\Gamma$ are contained in $N$.\\

By extending the techniques of the above theorem to asymptotic Plateau problem in $\BH^3$, we got the
following result.\\

\noindent \textbf{Theorem 4.5.} Let $\Gamma\subset \Si$ be a simple closed curve. Then either there exist a unique
minimal plane $\Sigma$ with $\PI \Sigma = \Gamma$, or there exist a canonical mean convex hull $N \subset \BH^3$ such
that all minimal planes $\Sigma '$ with $\PI \Sigma ' = \Gamma$ are contained in $N$.\\

By applying the above theorem to the limit set of a quasi-Fuchsian hyperbolic $3$-manifold, we got the
following corollary.\\

\noindent \textbf{Corollary 4.6.} Let $M$ be a quasi-Fuchsian hyperbolic $3$-manifold. Then either there exist a unique
minimal surface homotopy equivalent to $M$ or there exist a canonical mean convex core $N\subset M$ containing
all minimal surfaces homotopy equivalent to $M$.\\

Finally, we have a generic uniqueness result for extreme curves in $\BR^3$.\\

\noindent \textbf{Theorem 5.10.} Let $\Omega$ be a $C^3$-smooth mean convex domain in $\BR^3$, and $A=
\{\alpha \in C^3(S^1,\partial \Omega) \ | \ \alpha$ embedding$\}$. Then there exist an open dense subset
$A'\subset A$ in $C^3$ topology, such that for any $\Gamma\in A'$, there exist a unique least area
disk with boundary $\Gamma$.\\

The organization of the paper is as follows. In Section 2, we will give basic definitions and results which
will be used throughout the paper. In Section 3, canonical mean convex hull results will be proved for
extreme curves. In section 4, we will extend these to the asymptotic curves for $\BH^3$. In Section 5, we
will prove a generic uniqueness result for the extreme curves in $\BR^3$. Finally, we will have some remarks
on the results in Section 6.

\subsection{Acknowledgements:}

I would like to thank Bill Meeks, Yair Minsky, and Peter Li for very helpful conversations.

\section{Preliminaries}

In this section, we will overview the basic definitions and results which we use in the following sections.

\begin{defn}A {\em minimal disk (plane)} is a disk (plane) such that the mean curvature is $0$ at every point.
A {\em least area disk} is a disk which has the smallest area among the disks with same boundary. A least
area disk is minimal, but the converse is not true in general as minimal disks are just "locally" area
minimizing. A {\em least area plane} is a plane such that any subdisk in the plane is a least area disk.
\end{defn}

\begin{defn} Let $\Omega$ be a compact submanifold of a Riemannian $3$-manifold $M$. Then $\Omega$ is a {\em mean convex
domain} if the following conditions hold.

\begin{itemize}

\item $\partial \Omega$ is piecewise smooth.

\item Each smooth subsurface of $\partial \Omega$ has nonnegative curvature with respect to inward normal.

\item In the singular curves of $\partial \Omega$, the angle between the neighboring surfaces is less than
$\pi$ (inward direction).

\end{itemize}

\end{defn}

\begin{defn}
$\Gamma\subset M$ is an {\em extreme curve} if it is a curve in the boundary of a mean convex domain in $M$.
\end{defn}

\begin{rmk}
The above definition for extreme curves are different from the usual definition. We will abuse the name for
our purposes. This definition is more general than the definition in the literature, which says an extreme
curve is a curve in the boundary of a {\em \underline{convex domain}} .
\end{rmk}

\begin{defn} The sequence $\{\Sigma_{i}\}$ of embedded surfaces in a Riemannian manifold $M$ converges to the surface $\Sigma $ if

\begin{itemize}

\item $\Sigma $ contains all the limit points of the sequence $\{\Sigma_i\}$, i.e.\\
 $\Sigma = \{x= \lim_{i\to \infty }x_{i} \mid x_{i}\in S_{i}$ and $\{x_{i}\}$ convergent in $M\}$

\item Given $x\in\Sigma$ and $x_i\in \Sigma_i$ as above. Then there exist embeddings $f:D^2\rightarrow
\Sigma$, and $f_i: D^2 \rightarrow \Sigma_i$ with $f(0)=x$, and $f_i(0)=x_i$ such that $f_i$ converges to $f$
in $C^\infty$ topology.

\end{itemize}

\end{defn}

Now, we can state the main facts which we use in the first part.

\begin{lem}$[MY1]$
Let $\Omega$ be a mean convex domain, and $\Gamma\subset\partial \Omega$ be a simple closed curve. Then there
exist a least area disk $\Sigma\subset \Omega$ with $\partial \Sigma = \Gamma$. Moreover, all such disks are
properly embedded in $\Omega$ and they are pairwise disjoint. Moreover, If $\Gamma_1, \Gamma_2 \subset
\partial \Omega$ are disjoint simple closed curves, then the least area disks $\Sigma_1, \Sigma_2$ spanning
$\Gamma_1, \Gamma_2$ are also disjoint.
\end{lem}

\begin{lem}$[HS]$
Let $\Omega$ be a mean convex domain and let $\{\Sigma_i\}$ be a sequence of embedded least area disks in $\Omega$.
Then there is a subsequence $\{\Sigma_j\}$ of $\{\Sigma_i\}$ such that $\Sigma_j \rightarrow \Sigma$, a countable
collection of embedded least area disks in $\Omega$.
\end{lem}

\section{Mean Convex Hulls}

In this section, we will show that for any simple closed extreme curve in a Riemannian $3$-manifold, there
exist a canonical neighborhood which is mean convex or there exist a unique least area disk spanning the
curve. The idea is simple. Let $\Gamma\subset \partial \Omega$ be an extreme curve, and
$\Gamma_i^\pm\subset\partial \Omega$ be two sequences converging to $\Gamma$ from different sides. Then the
induced sequences of least area disks $\Sigma_i^\pm\subset\Omega$ with $\partial\Sigma_i^\pm = \Gamma_i^\pm$
limits to the two least area disks $\Sigma^+$ and $\Sigma^-$ with $\Sigma^\pm = \Gamma$. These least area
disks will be a barrier for the other least area disks with same boundary, and they will define a canonical
neighborhood $N$ with $\partial N = \Sigma^+ \cup \Gamma \cup \Sigma^-$, which we will call the {\em mean
convex hull} of $\Gamma$.

The combination of the above two lemmas gives us the following lemma, which is the main tool of this part.

\begin{lem}
Let $\Gamma\subset\partial \Omega$ be a simple closed extreme curve where $\Omega$ is the mean convex domain
in a $3$-manifold $M$. Then there are uniquely defined two canonical extremal least area disks (which might
be same) $\Sigma^+$ and $\Sigma^-$ in $\Omega$ with $\partial \Sigma^\pm = \Gamma$, and they are limits of
sequences of least area disks.

\end{lem}

\begin{pf}
$\Gamma\subset\partial \Omega$ is a simple closed extreme curve, and $\Omega$ is the mean convex domain in a
$3$-manifold $M$. Take a small neighborhood $A$ of $\Gamma$ in $\partial \Omega$, which will be a thin annulus with
$\Gamma$ is the core. $\Gamma$ separates the annulus $A$ into two parts, say $A^+$ and $A^-$ by giving a local
orientation. Define a sequence of pairwise disjoint simple closed curves $\Gamma_i^+\subset A^+ \subset \partial
\Omega$ such that $\lim\Gamma_i^+ = \Gamma$. Now, by Lemma 2.1, for any curve $\Gamma_i^+$, there exist an embedded
least area disk $\Sigma_i^+$ with $\partial \Sigma_i^+ = \Gamma_i^+$. This defines a sequence of least area disks
$\{\Sigma_i^+\}$ in $\Omega$. By Lemma 2.2, there exist a subsequence $\{\Sigma_j^+\}$ converging to a countable
collection of least area disks $\widehat{\Sigma}^+$ with $\partial \widehat{\Sigma}^+ = \Gamma$.

We claim that this collection $\widehat{\Sigma}^+$ consists of only one least area disk. Assume that there are two
disks in the collection say $\Sigma_a^+$ and $\Sigma_b^+$, and say $\Sigma_a^+$ is "above" $\Sigma_b^+$. By Lemma 2.1,
$\Sigma_a^+$ and $\Sigma_b^+$ are embedded and disjoint. They have same boundary $\Gamma\subset \Omega$. $\Sigma_b^+$
is also limit of the sequence $\{\Sigma_i^+\}$. But, since for any least area disk $\Sigma_i^+\subset \Omega$,
$\partial \Sigma_i^+ =\Gamma_i^+$ is disjoint from $\partial \Sigma_a^+ = \Gamma$, $\Sigma_i^+$ disjoint from
$\Sigma_a^+$, by exchange roundoff trick. This means $\Sigma_a^+$ is a barrier between the sequence $\{\Sigma_i^+\}$
and $\Sigma_b^+$, and so, $\Sigma_b^+$ cannot be limit of this sequence. This is a contradiction. So
$\widehat{\Sigma}^+$ is just one least area disk, say $\Sigma^+$. Similarly, $\widehat{\Sigma}^- = \Sigma^-$.

Now, we claim these least area disks $\Sigma^+$ and $\Sigma^-$ are canonical, depending only on $\Gamma$ and
$\Omega$, and independent of the choice of the sequence $\{\Gamma_i\}$ and $\{\Sigma_i\}$. Assume that there
exist another least area disk $S^+$ with $\partial S^+ = \Gamma$ and $S^+$ is a limit of the sequence of
least area disks $S_i^+$ with $\partial S_i^+ = \gamma_i^+ \subset A^+$. By Lemma 2.1, $\Sigma^+$ and $S^+$
are disjoint. Then one of them is "above" the other one. If $\Sigma^+$ is above $S^+$, then $\Sigma^+$
between the sequence $S_i^+$ and $S^+$. This is because, all $S_i^+$ are disjoint and above $S^+$ as
$\partial S_i^+ = \gamma_i$ are disjoint and "above" $\Gamma$. Similarly, $\Sigma^+$ is below $S_i$ for any
$i$, as $\partial \Sigma^+ = \Gamma$ is below the curves $\gamma_i^+\subset A^+$. Now, since $\Sigma^+$ is
between the sequence $\{S_i^+\}$ and its limit $S^+$, and $S^+$ and $\Sigma^+$ are disjoint, $\Sigma^+$ will
be a barrier for the sequence $\{S_i^+\}$, and so it cannot limit on $S^+$. But, this is a contradiction.
Similarly, $\Sigma^+$ cannot be below $S^+$, so they must be same. So, $\Sigma^+$ and $\Sigma^-$ are
canonical least area disks for $\Gamma$.
\end{pf}

\begin{thm}
Let $\Omega$ be a mean convex domain in a $3$-manifold $M$ and $\Gamma\subset\partial \Omega$ be a simple
closed curve. Then either there exist a unique least area disk $\Sigma$ in $\Omega$ with $\partial \Sigma =
\Gamma$, or there exist a canonical mean convex hull $N$  in $\Omega$ such that $\partial N= \Sigma^+ \cup
\Sigma^-$ where $\Sigma^\pm$ are uniquely defined extremal least area disks in $\Omega$ with $\partial
\Sigma^\pm = \Gamma$. Moreover, all least area disks $\Sigma '\subset \Omega$ spanning $\Gamma$ are contained
in $N$.
\end{thm}

\begin{pf}
$\Gamma\subset\partial \Omega$ is a simple closed extreme curve, and $\Omega$ is the mean convex domain in a
$3$-manifold $M$. Let $\Sigma^+$ and $\Sigma^-$ be the extremal least area disks for $\Gamma$ from Lemma 3.1.
Let $N \subset M$ be the region between $\Sigma^+$ and $\Sigma^-$, i.e. $\partial N= \Sigma^+ \cup \Sigma^-$.
Assume $\Sigma'$ is a least area disk with $\partial \Sigma' = \Gamma$. We claim that $\Sigma' \subset N$.
Assume on the contrary. Since all least area disks are disjoint by Lemma 2.1, $\Sigma'\cap\Sigma^\pm =
\emptyset$, which implies $\Sigma'\cap N =\emptyset$. Then either $\Sigma'$ is "above" $\Sigma^+$ or "below"
$\Sigma^-$. If $\Sigma'$ is "above" $\Sigma^+$, then since there is a sequence $\{\Sigma_i^+\}$ such that
$\Sigma_i^+ \rightarrow \Sigma^+$, for sufficiently large $k$, $\Sigma_k^+ \cap \Sigma' \neq \emptyset$ as
$\partial \Sigma_k^+= \Gamma_k^+$ is above $\partial \Sigma = \Gamma$. But by the choice of the sequence,
$\Gamma_k^+\cap \Gamma =\emptyset$. So, by Lemma 2.1, the least area disks $\Sigma_k^+$ and $\Sigma'$ must be
disjoint. This is a contradiction. Similarly, $\Sigma'$ cannot be "below" $\Sigma^-$. So $\Sigma'\subset N$.

If $\Sigma^+ = \Sigma^-$, then $N = \Sigma^+ = \Sigma^-$. Since for any least area disk $\Sigma'\subset
\Omega$ with $\partial \Sigma' = \Gamma$, $\Sigma'$ is contained in $N$, then $\Sigma' = \Sigma^+ =
\Sigma^-$. This means there exist a unique least area disk spanning $\Gamma$.
\end{pf}

\begin{rmk}
In Lemma 3.1, and Theorem 3.2, the canonical least area disks, and mean convex hull for $\Gamma \subset
\partial\Omega$ are depending also on the mean convex domain $\Omega$. This is because, even though $\Gamma$
is in the boundary of mean convex domain $\Omega$, it does not mean that any least area disk in $M$ spanning
$\Gamma$ must be in $\Omega$.

We call the region $N$ assigned to $\Gamma$ as {\em mean convex hull}. The reason for this $N\subset\Omega$
is itself a mean convex domain. The importance of this object is that it is canonically defined and uniquely
determined by $\Gamma$ and $\Omega$. One can think of this object as a {pseudo-convex hull} living in
$\Omega$ and in the convex hull of $\Gamma$.

If you want the canonical least area disks $\Sigma^\pm$, and mean convex hull $N$ for $\Gamma$ to be
independent from the mean convex domain $\Omega$, you need to make sure that all least area disks spanning
$\Gamma$ are in the mean convex domain. One condition to guarantee that is that $\Gamma$ has convex hull
property and the convex hull of $\Gamma$ is in $\Omega$. If $\Omega$ is convex domain, this is automatic. In
other words, If $\Omega$ is convex domain (i.e. $\Gamma$ is an extreme curve in the usual sense), than the
defined mean convex hull, and extremal least area disks are independent of $\Omega$, depends only on
$\Gamma$. Moreover, all least area disks in $M$ spanning $\Gamma$ are in the mean convex hull, i.e. between
the extremal least area disks.

Also, a similar result was obtained by Brian White by using geometric measure theory methods in [Wh].
\end{rmk}

\section{Mean Convex Hulls in Hyperbolic Space}

Now, we are going to generalize the above results to the least area planes in hyperbolic space. Our aim in this section
to show the existence of canonical mean convex hulls for a simple closed curve in $\Si$. The technique is basically
same. To prove this, we need analogies of the lemmas in Section 2.

\begin{lem}$[An2]$
Let $\Gamma\subset \Si$ be a simple closed curve in the sphere at infinity of hyperbolic $3$-space. Then there exist a
properly embedded least area plane $\Sigma$ spanning $\Gamma$, i.e. $\PI \Sigma = \Gamma$.
\end{lem}

\begin{lem}$[Ga]$
Let $\{\Sigma_i\}$ be a sequence of least area planes in $\BH^3$ with $\PI \Sigma_i = \Gamma_i \subset \Si$ simple
closed curve for any $i$. If $\Gamma_i\rightarrow\Gamma$, then there exist a subsequence $\{\Sigma_j\}$ of
$\{\Sigma_i\}$ such that $\Sigma_j\rightarrow\Sigma$ a collection of least area planes whose asymptotic boundaries are
$\Gamma$.
\end{lem}

\begin{lem}
Let $\Gamma_1$ and $\Gamma_2$ be two disjoint simple closed curves in $\Si$. Then, if $\Sigma_1$ and $\Sigma_2$ are
least area planes with $\PI \Sigma_i = \Gamma_i$, then $\Sigma_1$ and $\Sigma_2$ are disjoint, too.
\end{lem}

\begin{pf}
Assume that $\Sigma_1\cap\Sigma_2\neq\emptyset$. Since asymptotic boundaries $\Gamma_1$, and $\Gamma_2$ are
disjoint, the intersection cannot contain an infinite line. So, the intersection between $\Sigma_1$ and
$\Sigma_2$ must contain a simple closed curve $\gamma$. Now, $\gamma$ bounds two least area disks $D_1$ and
$D_2$ in $\BH^3$, with $D_i\subset\Sigma_i$. Now, take a larger subdisk $E_1$ of $\Sigma_1$ containing $D_1$,
i.e. $D_1\subset E_1 \subset \Sigma_1$. By definition, $E_1$ is also a least area disk. Now, modify $E_1$ by
swaping the disks $D_1$ and $D_2$. Then, we get a new disk $E_1 '= \{E_1 - D_1\} \cup D_2$. Now, $E_1$ and
$E_1 '$ have same area, but $E_1 '$ have folding curve along $\gamma$. By smoothing out this curve as in
[MY1], we get a disk with smaller area, which contradicts to $E_1$ being least area. Note that this technique
is known as Meeks-Yau exchange roundoff trick.
\end{pf}

Now, we will adapt Lemma 3.1 to our context.

\begin{lem}
Let $\Gamma\subset \Si$ be a simple closed curve. Then there are two canonical extremal least area planes
(which might be same) $\Sigma^+$ and $\Sigma^-$ in $\BH^3$ with $\PI \Sigma^\pm = \Gamma$. Moreover, any
least area plane $\Sigma$ with $\PI\Sigma = \Gamma$ is disjoint from $\Sigma^\pm$, and it is captured in the
region bounded by $\Sigma^+$ and $\Sigma^-$ in $\BH^3$.
\end{lem}

\begin{pf}
Let $\Gamma\subset \Si$ be a simple closed curve. $\Gamma$ separates $\Si$ into two parts, say $D^+$ and $D^-$. Define
sequences of pairwise disjoint simple closed curves $\{\Gamma_i^+\}$ and $\{\Gamma_i^-\}$ such that $\Gamma_i^+\subset
D^+$, and $\Gamma_i^-\subset D^-$ for any $i$, and $\Gamma_i^+ \rightarrow \Gamma$, and $\Gamma_i^- \rightarrow
\Gamma$.

By Lemma 4.1, for any $\Gamma_i^+\subset \Si$, there exist a least area plane $\Sigma_i^+\subset \BH^3$. This defines a
sequence of least area planes $\{\Sigma_i^+\}$. Now, by using Lemma 4.2, we get a collection of least area planes
$\widehat{\Sigma}^+$ with $\PI\widehat{\Sigma}^+ = \Gamma$, as $\PI\Sigma_i^+ = \Gamma_i^+ \rightarrow \Gamma$.

As in the proof of Lemma 3.1, we claim that the collection $\widehat{\Sigma}^+$ consists of only one least
area plane. But this time, we do not know that if two least area planes have same asymptotic boundary, then
they they are disjoint, as in the compact domain case (Lemma 2.1). We only know that if two least area planes
have disjoint asymptotic boundary, then they are disjoint by Lemma 4.3.

Assume that there are two least area planes $\Sigma_a^+$ and $\Sigma_b^+$ in the collection
$\widehat{\Sigma}^+$. Since $\PI \Sigma_a^+ = \PI \Sigma_b^+ = \Gamma$, $\Sigma_a^+$ and $\Sigma_b^+$ might
not be disjoint, but they are disjoint from least area planes in the sequence, i.e. $\Sigma_i^+\cap
\Sigma_{a,b}^+ =\emptyset$ for any $i$, by Lemma 4.3.

If $\Sigma_a^+$ and $\Sigma_b^+$ are disjoint, we can use the argument in the proof Lemma 3.1, and conclude
that $\Sigma_a^+$ is a barrier between the sequence $\{\Sigma_i^+\}$ and $\Sigma_b^+$, and so, $\Sigma_b^+$
cannot be limit of this sequence.

If $\Sigma_a^+$ and $\Sigma_b^+$ are not disjoint, then they intersect each other, and in some region,
$\Sigma_b^+$ is "above" $\Sigma_a^+$. But since $\Sigma_a^+$ is the limit of the sequence $\{\Sigma_i^+\}$,
this would imply $\Sigma_b^+$ must intersect planes $\Sigma_i^+$ for sufficiently large $i$. But, this
contradicts $\Sigma_b^+$ is disjoint from $\Sigma_i^+$ for any $i$, as they have disjoint asymptotic
boundary. So, there exist unique least area plane $\Sigma^+$ in the collection $\widehat{\Sigma}^+$.
Similarly, $\widehat{\Sigma}^- = \Sigma^-$.

By using similar arguments to Lemma 3.1, one can conclude that these least area planes $\Sigma^+$, and
$\Sigma^-$ are canonical, and independent of the choice of the sequence $\{\Gamma_i^\pm\}$ and
$\{\Sigma_i^\pm\}$.

Now, if we show the last statement of the theorem, then we are done. Let $\Sigma'$ be any least area plane
with $\PI\Sigma'=\Gamma$. If $\Sigma'\cap\Sigma^+ \neq \emptyset$, then some part of $\Sigma'$ must be
"above" $\Sigma^+$. Since $\Sigma^+=\lim \Sigma_i^+$, for sufficiently large $i$, $\Sigma'\cap\Sigma_i^+
\neq\emptyset$. But, $\PI\Sigma_i^+ = \Gamma_i^+$ is disjoint from $\Gamma=\PI\Sigma'$. Then, by Lemma 4.3,
$\Sigma'$ must be disjoint from $\Sigma_i^+$, which is a contradiction.

Similarly, this is true for $\Sigma^-$, too. Moreover, let $N \subset \BH^3$ be the region between $\Sigma^+$
and $\Sigma^-$, i.e. $\partial N= \Sigma^+ \cup \Sigma^-$. Then, by construction,  $N$ is also a canonical
region for $\Gamma$, and for any least area plane $\Sigma'$ with $\PI \Sigma'=\Gamma$, $\Sigma'$ is contained
in the region $N$, i.e. $\Sigma\subset N$.
\end{pf}

\begin{rmk}
This lemma shows that for any given simple closed curve $\Gamma$ in $\Si$, we can get two canonical least
area planes $\Sigma_\Gamma^\pm$ (which might be same). Moreover, these least area planes are disjoint from
any other least area plane with asymptotic boundary $\Gamma$. Because of this reason, we call these canonical
least area planes $\Sigma_\Gamma^\pm$ as {\em untouchable least area planes}. Also, the region between them
is a canonical mean convex hull of $\Gamma$, which captures all least area planes in $\BH^3$ spanning
$\Gamma$.
\end{rmk}

\begin{thm}
Let $\Gamma\subset \Si$ be a simple closed curve. Then either there exist a unique minimal plane $\Sigma$ with $\PI
\Sigma = \Gamma$, or there exist a canonical mean convex hull $N \subset \BH^3$ such that all minimal planes $\Sigma '$
with $\PI \Sigma ' = \Gamma$ are contained in $N$.
\end{thm}

\begin{pf}
Let $\Gamma\subset \Si$ be a simple closed curve and let $X_\Gamma =\{P\subset\BH^3 \ | \ \PI P = \Gamma, P
\mbox{ is a minimal plane} \}$. By [An1], we know that for any $P\in X_\Gamma$, $P$ is in the convex hull of
its asymptotic boundary, i.e. $P\subset CH(\Gamma)$. Also, name the complement of $\Gamma$ in $\Si$ as
$D^\pm$, i.e. $\Si- \Gamma = D^+ \cup D^-$.

Now, for any $P_i\in X_\Gamma$, define domains $\Delta_i^\pm\subset \BH^3$ such that $\Delta_i^+$ is the
unbounded component of $\BH^3 - P_i$ with asymptotic boundary $D^+\subset \Si$. $\Delta^-$ is defined
similarly. Now, let $\Delta^+=\bigcap_i \Delta_i^+$, and $\Delta^-$ is defined similarly. Now, $\Delta^\pm$
are nonempty as all minimal planes stays in convex hull of $\Gamma$. Also, $\Delta^\pm$ are canonical and
mean convex as boundaries coming from minimal planes.

Now, as in Lemma 4.4, we will define extremal least area planes in $\Delta^+$. Now, the crucial point here is
that they are no more least area planes in $\BH^3$, but least area planes in $\Delta^+$.

As in the proof of Lemma 4.4, we will start with two sequences of simple closed curves converging to $\Gamma$
from different sides. Define sequences of pairwise disjoint simple closed curves $\{\Gamma_i^+\}$ and
$\{\Gamma_i^-\}$ such that $\Gamma_i^+\subset D^+$, and $\Gamma_i^-\subset D^-$ for any $i$, and $\Gamma_i^+
\rightarrow \Gamma$, and $\Gamma_i^- \rightarrow \Gamma$.

Since, $\Delta^+$ is mean convex domain, we can define a least area plane $\Sigma_i^+$ in $\Delta^+$ with
asymptotic boundary $\Gamma_i^+$ as follows. Take a sequence of extremal simple closed curves in $\Delta^+$
such that these curves converge to $\Gamma_i^+$. Then by [MY2], there exist least area disks in the mean
convex domain $\Delta^+$ spanning the extremal simple closed curves in the sequence. Then by taking a limit
of the sequence of these least area disks, we will get a least area plane in $\Delta^+$ as in [Ga]. Note that
this least area plane is not a least area plane of $\BH^3$, but least area plane of $\Delta^+$.

Now, as in Lemma 4.4, we will take the limit of least area planes $\Sigma_i^+$, and we get uniquely defined
limit least area plane $\Sigma^+$ in $\Delta^+$ by Lemma 4.4. Similarly, construct the least area plane
$\Sigma^-$ in $\Delta^-$. Note that even though $\Sigma^\pm$ are canonical least area planes in $\Delta^\pm$,
they might not have least area property in the whole $\BH^3$.

Now, let $\Sigma^+$ and $\Sigma^-$ be the canonical planes for $\Gamma$ as above, and let $N \subset \BH^3$
be the region between $\Sigma^+$ and $\Sigma^-$, i.e. $\partial N= \Sigma^+ \cup \Sigma^-$. We claim that for
any minimal plane $P_i\in X_\Gamma$, $P_i\subset N$. Indeed, this is clear by construction. $\Sigma^+\subset
\Delta^+$ and by definition $\Delta^+=\bigcap_i \Delta_i^+$. Since $P_i$ is "below" $\Delta_i^+$, then it is
"below" $\Delta^+$, and so it is "below" $\Sigma^+$. Similarly $P_i$ is "above" $\Delta^-$ and so it is above
$\Sigma^-$. This implies that $P_i\subset N$.

If $\Sigma^+ = \Sigma^-$, then $N = \Sigma^+ = \Sigma^-$. For any minimal plane with $\PI P = \Gamma$,
$\Sigma$ is contained in $N$. So, $P = \Sigma^+ = \Sigma^-$. This means there exist a unique minimal plane
$P\subset\BH^3$ spanning $\Gamma\subset \Si$. Moreover, this unique minimal plane is also least area by
existence of least area planes by [An2].

\end{pf}

\begin{rmk}
This result is also true for cocompact Gromov hyperbolic spaces. If $M$ is a compact Gromov hyperbolic $3$-manifold,
then the above statements are true for the universal cover $\widetilde{M}$. In other words, if you replace $\Si$ with
$S^2_\infty(\widetilde{M})$ and $\BH^3$ with $\widetilde{M}$ in the statement of Theorem 4.5, it is still true. This is
because one can prove Lemma 4.4 by modifying the results in [Co1], and the proof of Theorem 4.5 simply goes through.
\end{rmk}

\subsection{Quasi-Fuchsian Manifolds}

\begin{defn}
Let $M$ be a hyperbolic $3$-manifold with $\pi_1 (M) = \pi_1 (S)$ for some compact oriented surface $S$. Then
$M$ is called {\em quasi-Fuchsian} if the limit set (asymptotic limit of the orbit of a point in
$\widetilde{M} = \BH^3$ under covering transformations) is a simple closed curve in $\Si$.
\end{defn}

\begin{cor} Let $M$ be a quasi-Fuchsian hyperbolic $3$-manifold. Then either there exist a unique minimal surface
homotopy equivalent to $M$ or there exist a canonical mean convex core $N\subset M$ containing all minimal
surfaces homotopy equivalent to $M$.
\end{cor}

\begin{pf}
We will use the notation of Theorem 4.5. Since $M$ is quasi-Fuchsian, its limit set is a simple closed curve,
say $\Gamma\subset\Si$. By Theorem 4.5, we have a canonical mean convex hull $\widehat{N}\subset\BH^3$ with
$\PI \widehat{N} =\Gamma$. Since $\widehat{N}$ is only depending on $\Gamma$, in order to get the canonical
mean convex core, all we need to show that $\widehat{N}$ is invariant under covering transformation. So, if
we can show that the planes $\partial \widehat{N} = \Sigma^\pm$ are invariant under covering transformations,
then we are done.

Let $G \simeq \pi_1(M)$ be the covering transformations of $M$. Then for any $\alpha\in G$, $\alpha$ induces
a homeomorphism on $\Si$ and $\alpha(\Gamma) = \Gamma$. Since $\Sigma^+$ is the limit of the sequence
$\{\Sigma_i^+\}$ with $\PI \Sigma_i^+ =\Gamma_i^+\subset D^+$, then $\alpha(\Sigma^+)$ is the limit of the
sequence $\{\alpha(\Sigma_i^+)\}$. Since $M$ is orientable, $\alpha(D^\pm)=D^\pm$ and $\alpha(\Gamma_i^+)
\subset D^+$. So $\{\Gamma_i^+\}$ is a sequence of pairwise disjoint simple closed curves in $D^+$, and
$\alpha(\Gamma_i^+) \rightarrow \Gamma$. But the proof of Lemma 4.4 implies there is only one limit for such
a sequence of least area planes, i.e. $\alpha(\Sigma^+) = \Sigma^+$. Similarly, $\alpha(\Sigma^-) =
\Sigma^-$. So $\widehat{N}$ is invariant under covering transformations. So, under covering projection $\pi :
\BH^3 \rightarrow M$, $\pi(\widehat{N})= N \subset M$ defines the desired canonical mean convex core.

Now, we claim that if $S\subset M$ is a minimal surface homotopy equivalent to $M$, then $S$ is contained in
our canonical mean convex core $N$. Since $S$ is homotopy equivalent to $M$, it is $\pi_1$-injective surface
in $M$ and its universal cover $\widetilde{F}\subset \BH^3$ is a minimal plane such that $\PI \widetilde{F} =
\Gamma$. Then, Theorem 4.5 implies that $\widetilde{F}\subset \widehat{N}$. So if we take the projection of
both, we get $F\subset N$.

In the case, $\widehat{N}=\Sigma^+ = \Sigma^-$, $N$ will be a least area surface, which is homotopy
equivalent to $M$. Since any minimal surface homotopy equivalent to $M$ contained in $N$, this implies there
exist a unique minimal surface homotopy equivalent to $M$. Moreover, this unique minimal surface is indeed
least area by existence of least area planes by [An2].
\end{pf}

\begin{rmk}
The {\em mean convex core} lives in the convex core of $M$. This is because all minimal planes $\Sigma$ with
$\PI \Sigma = \Gamma$, $\Sigma$ is contained in the convex hull of $\Gamma$, say $CH(\Gamma)$ (smallest
convex subset of $\BH^3$ with asymptotic boundary $\Gamma\subset \Si$). Then the mean convex hull is in the
convex hull, $\widehat{N} \subset CH(\Gamma)$, which implies mean convex core is in the convex core of $M$,
since the convex core of $M$ is the projection of $CH(\Gamma)$.
\end{rmk}

\section{Generic Uniqueness}

In this section, we will give a generic uniqueness result for least area disks spanning an extreme curve. In
other words, we will show that a generic curve on the boundary of a mean convex domain bounds a unique least
area disk. Similar result for least area planes in $\BH^3$ has been proved in [Co3].

Let $\Omega$ be a $C^3$-smooth mean convex domain in $\BR^3$. We need to define the following spaces.

$A= \{\alpha \in C^3(S^1,\partial \Omega)\ | \ \  \alpha$ embedding $\}$

$D=\{u \in C^3(S^1,S^1)\ |\ \  u$ diffeomorphism and satisfies three point condition, i.e. $u(e^{\frac{2}{3}
k \pi i}) = e^{\frac{2}{3} k \pi i}, k=1,2,3 \}$

$M=\{f:D^2\rightarrow \BR^3 \  | \ \  f(D^2)$ minimal and $ f|_{\partial D^2} \in A \}$

Now, we will quote Tomi and Tromba's results from [TT] on the structure of these spaces. They consider a
minimal map $f:D^2\rightarrow \BR^3$, as conformal harmonic map, and realize the space of minimal maps as a
subspace of space of harmonic maps. The minimal ones correspond to the conformal ones in this space. On the
other hand, one can identify a harmonic map from a disk to $\BR^3$ with its boundary parametrization, by
unique extension property. So, we can think of the space of minimal maps $M$, as a subspace of harmonic maps
or their boundary parametrizations, $A$ in above notation. If you don't care about the parametrization but
the image curve, you can augment the space of boundary parametrizations with a "reparametrization" factor $D$
to capture conformality.

So, they consider $A \times D$ such that $(\alpha, u )\in A\times D$ identified with $\widetilde {\alpha
\circ u} :D^2\rightarrow \BR^3$, the harmonic extension of $\alpha \circ u: S^1\rightarrow \BR^3$. Then, any
minimal disk spanning $\alpha(S^1) \subset \partial \Omega$, will correspond to a point in the slice
$\{\alpha\} \times D$ corresponding to a {\em conformal} harmonic map. They define a conformality operator $k
: A \times D \rightarrow C^2(S^1)$ such that the kernel of this map will be the conformal maps. So, by
abusing the above notation, define $M := ker(k) \subset A \times D$, the space of minimal disks with boundary
in $\partial \Omega$.

In [TT], Tomi and Tromba proved that the second component of the derivative of conformality operator, $D_u k:
T_u D \rightarrow Z\subset C^2(S^1)$ is almost an isomorphism, which is a Fredholm map of index $0$. Then by
using basic linear algebra they show that the projection maps $\Pi_1 :A\times D \rightarrow A$ restricted to
the $M=ker(k)$, i.e. $\Pi_1 |_M: M\rightarrow A$ is Fredholm of index $0$. This means the restriction map
from minimal maps to their boundary parametrizations is almost an isomorphism.

Now, in our case, $A$ is different from Tomi and Tromba's setting. But, if one look at the Tomi and Tromba's
proof, everything happens in the second component of $A\times D$, and they get the "almost isomorphism"
between $D$ and image of $k$. So, by simple modifications of these proofs would give the following desired
result for our purposes. To see the modifications in detail, one can look at [Co2], and [Co3].

\begin{lem} $[TT]$
Let $A, D$, and $M$ be as above, and $\Pi_1:A\times D \rightarrow A$ be the projection map. Then the space of
minimal maps $M$, is a submanifold of $A\times D$, and restriction of projection map to $M$, $\Pi_1 |_M$ is
Fredholm of index $0$.
\end{lem}

The second lemma which we use is the classical inverse function theorem for Banach manifolds, [La].

\begin{lem}
(Inverse Function Theorem) Let $M$ and $N$ be Banach manifolds, and let $F: M \rightarrow N$ be a $C^p$ map.
Let $x_0 \in M$ and $dF$ is isomorphism at $x_0$. Then $F$ is local $C^p$ diffeomorphism, i.e. there exist an
open neighborhood of $U\subset M$ of $x_0$ and an open neighborhood $V \subset N$ of $F(x_0)$ such that
$F|_U:U\rightarrow V$ is $C^p$ diffeomorphism.
\end{lem}

The last ingredient is the generalization of Sard's theorem to infinite dimensional spaces [Sm].

\begin{lem}(Sard-Smale Theorem)
Let $F: X\rightarrow Y$ be a Fredholm map. Then the regular values of $F$ are almost all of $Y$, i.e except a
set of first category.
\end{lem}

From now on, we will call the regular values of the Fredholm map $\Pi_1 |_M$ as generic curve. Now, we can
establish the main analytical tool for our generic uniqueness result.

\begin{lem}
Let $\alpha\in A$ be a generic curve. Then for any $\Sigma\in \Pi_1^{-1}(\alpha)$, there exist a neighborhood
$U_\Sigma\subset M$ such that $\Pi_1 |_{U_\Sigma}$ is a homeomorphism onto a neighborhood of $\alpha$ in $A$.
\end{lem}

\begin{pf} By Lemma 5.1, the map $\Pi_1 |_M: M \rightarrow A$ is Fredholm of index 0. Let $\alpha\in A$ be a generic curve
and $\Sigma \in \Pi_1^{-1}(\alpha)\subset M$. Since $\alpha$ is regular value, $D\Pi_1(\Sigma):T_\Sigma M
\rightarrow T_\alpha A$ is surjective, i.e. $dim(\mbox{coker}(D\Pi_1))= 0$. Moreover, we know that $\Pi_1$ is
Fredholm of index 0. This implies $dim(\mbox{ker}(D\Pi_1))=dim(\mbox{coker}(D\Pi_1))= 0$, and so $D\Pi_1$ is
isomorphism at the point $\Sigma\in M$. By the Inverse Function Theorem, there exist a neighborhood of
$\Sigma$ which $\Pi_1$ maps homeomorphically onto a neighborhood of $\alpha\in A$.
\end{pf}

Now, our aim is to construct a foliated neighborhood for any least area disk spanning a given generic curve
(regular value of the Fredholm map) in $A$. Moreover, we will show that the leaves of this foliation are
embedded least area disks with pairwise disjoint boundary. By using this, we will show that uniqueness of the
least area disk spanning the generic curve.

We will abuse the notation by using interchangeably the map $\Gamma: S^1\rightarrow \partial \Omega$ with its
image $\Gamma(S^1)$. Similarly same is true for $\Sigma :D^2\rightarrow \Omega$ and its image $\Sigma(D^2)$.

Let $\Gamma_0\in A$ be a generic curve, and let $\Sigma_0\in \pi_1^{-1}(\Gamma_0)\subset M$ be a least area
disk whose existence guaranteed by [MY1]. Then by Lemma 5.4, there is a neighborhood of $\Sigma_0\in U\subset
M$ homeomorphic to the neighborhood $\Gamma_0\in V \subset A$.

Let $\Gamma:[-\epsilon,\epsilon]\rightarrow V$ be a path such that $\Gamma(0)=\Gamma_0$ and for any $t, t'
\in[-\epsilon,\epsilon]$, $\Gamma_t\cap\Gamma_{t'}=\emptyset$. In other words, $\{\Gamma_t\}$ foliates a
neighborhood of $\Gamma_0$ in $\partial \Omega$. Let $\Sigma_t\in U$ be the preimage of $\Gamma_t$ under the
local homeomorphism.

\begin{lem}
$\{\Sigma_t\}$ is a foliation of a neighborhood of $\Sigma_0$ in $\Omega$ by embedded least area disks.
\end{lem}

\begin{pf} We will prove the lemma in three steps.

\vspace{0.3cm}

\textbf{Claim 1:} For any $s\in [-\epsilon, \epsilon]$, $\Sigma_s$ is an embedded disk.

\vspace{0.1cm}

\begin{pf}
Since $\Sigma_0$ is a least area disk, by Lemma 2.1, $\Sigma_0$ is an embedded disk. Now, $\{\Sigma_t\}$ is
continuous family of minimal disks. We cannot apply the Lemma 2.1 to these disks, since the lemma is true for
least area disks, while our disks are only minimal.

Let $s_0=inf\{s\in(0,\epsilon] \ | \ \Sigma_s \mbox{ is not embedded}\}$. But, since $\{\Sigma_t\}$ is
continuous family of disks, and this can only happen when $\Sigma_{s_0}$ has tangential self intersection
(locally lying on on side). But this contradicts to maximum principle for minimal surfaces. So for all $s\in
[0,\epsilon]$, $\Sigma_s$ is embedded. Similarly, this is true for $s\in [-\epsilon, 0]$, and the result
follows.
\end{pf}

\textbf{Claim 2:} $\{\Sigma_t\}$ is a foliation, i.e. for any $t, t' \in[\epsilon,\epsilon]$,
$\Sigma_t\cap\Sigma_{t'}=\emptyset$.

\vspace{0.1cm}

\begin{pf}
Assume on the contrary that there exist $t_1 < t_2$ such that $\Sigma_{t_1}\cap\Sigma_{t_2}\neq \emptyset$.
First, since the boundaries $\Gamma_{t_1}$ and $\Gamma_{t_2}$ are disjoint, the intersection cannot contain a
line segment. So the intersection must be a collection of closed curves. We will show that in this situation,
there must be a tangential intersection between two disks, and this will contradict to the maximum principle
for minimal surfaces.

If $\Sigma_{t_2}$ does not intersect all the minimal disks $\Sigma_s$ for $s\in [-\epsilon, t_2]$, let
$s_0=\sup\{s\in [-\epsilon, t_2] \ | \ \Sigma_{t_2}\cap\Sigma_s = \emptyset\}$. Then, since $\{\Sigma_t\}$ is
continuous family of minimal disks, it is clear that $\Sigma_{t_2}$ must intersect $\Sigma_{s_0}$
tangentially, and lie in one side of $\Sigma_{s_0}$. But this contradicts to maximum principle for minimal
surfaces.

So, let's assume $\Sigma_{t_2}$ intersects all minimal disks $\Sigma_s$ for $s\in [-\epsilon, t_2]$. Let
$s_0=\sup\{s\in[-\epsilon,+\epsilon] \ | \ \Sigma_{-\epsilon} \cap \Sigma_s =\emptyset\}$. If $s_0 >
-\epsilon$, then this would imply a tangential intersection as above, which is a contradiction. Otherwise,
the supremum is $-\epsilon$, which implies for any $t > -\epsilon$, $\Sigma_t$ intersects
$\Sigma_{-\epsilon}$. In particular, this implies $\Sigma_0\cap \Sigma_{-\epsilon} \neq \emptyset$. Now,
$\Sigma_{-\epsilon}$ separates $\Omega$, and defines a mean convex domain. But, $\Gamma_0$ is in the boundary
of this new mean convex domain, and the least area disk $\Sigma_0$ must be embedded inside of this mean
convex domain by Lemma 2.1. So, $\Sigma_0$ cannot intersect $\Sigma_{-\epsilon}$. This is a contradiction,
and the result follows.
\end{pf}

\textbf{Claim 3:} For any $s\in [-\epsilon, \epsilon]$, $\Sigma_s$ is a least area disk.

\vspace{0.1cm}

\begin{pf} Fix $\Sigma_s$ for $s\in (-\epsilon, \epsilon)$. Now, let $[\Sigma_{-\epsilon}, \Sigma_{\epsilon}]$ be the region
bounded by embedded disks $\Sigma_{-\epsilon}$ and $\Sigma_{\epsilon}$ in $\Omega$. By above results,
$\Sigma_s\subset [\Sigma_{-\epsilon}, \Sigma_{\epsilon}]$. Since the boundaries are minimal disks,
$[\Sigma_{-\epsilon}, \Sigma_{\epsilon}]$ is a mean convex region. Let $\gamma\subset \Sigma_s$ be a simple
closed curve. By [MY1], there exist a least area embedded disk $D$ spanning $\gamma$ in the mean convex
domain $[\Sigma_{-\epsilon}, \Sigma_{\epsilon}]$. If $D$ is not in $\Sigma_s$, it must intersect other leaves
nontrivially. Then $\{\Sigma_t\}\cap D$ induce a singular 1-dimensional foliation $F$ on $D$. The
singularities of the foliation are isolated as $\{\Sigma_t\}$ and $D$ are minimal disks. Since Euler
characteristic of the disk is 1, by Poincare-Hopf index formula there must be a positive index singularity
implying tangential (lying on one side) intersection of $D$ with some leave $\Sigma_s$. But this contradicts
to maximum principle for minimal surfaces. Since $\epsilon$ was chosen arbitrarily at the beginning, one can
start with suitable $\epsilon'>\epsilon$. The whole proof will go through, and this shows that for any $s\in
[-\epsilon, \epsilon]$, $\Sigma_s$ is a least area disk.
\end{pf}
\end{pf}

Existence of such a foliated neighborhood for a least area disk, implies uniqueness.

\begin{lem}
$\Sigma_0$ is the unique least area disk with boundary $\Gamma_0$.
\end{lem}

\begin{pf}  Let $\Sigma'$ be another least area disk with boundary $\Gamma_0$. If $\Sigma_0 \neq \Sigma'$ then
$\Sigma'$ must intersect a leave in the foliated neighborhood of $\Sigma_0$, say $\Sigma_s$. But, since
$\partial\Sigma_s = \Gamma_s$ is disjoint from $\partial\Sigma'=\Gamma_0$, Lemma 2.1 implies that the least area disks
$\Sigma_s$ and $\Sigma '$ are disjoint. This is a contradiction.
\end{pf}

So, we have proved the following theorem:

\begin{thm}
Let $\Gamma\in A$ be a generic curve as described above. Then there exist a unique least area disk
$\Sigma\subset \Omega$ with $\partial\Sigma=\Gamma$.
\end{thm}

\begin{rmk} This theorem does not say that there exist a unique {\it minimal} disk spanning a given generic curve.
In the proof of Lemma 5.5, we essentially use the disk $\Sigma_0$ being least area.
\end{rmk}

So far we have proved the uniqueness of least area disks for a subset $\widehat{A}\subset A$, where $A -
\widehat{A}$ is a set of first category. In the following subsection, we will show that this is true for a
more general class of curves, i.e. an open dense subset of $A$.

\subsection{Open dense set of curves}

Now, we will show that any regular curve has an open neighborhood such that the uniqueness result holds for
any curve in this neighborhood.

Let $\Gamma_0\in A$ be a regular curve, and let $\Sigma_0\in \pi_1^{-1}(\Gamma_0)\subset M$ be the unique
least area disk spanning $\Gamma_0$. Let $ U\subset M$ be the neighborhood of $\Sigma_0$ homeomorphic to the
neighborhood $V \subset A$ of $\Gamma_0$ as above. We will show that $\Gamma_0$ has a smaller open
neighborhood $V'\subset V$ such that for any $\Gamma\in V'$, there exist unique least area disk in $\Omega$
with $\partial\Sigma=\Gamma$.

First we will show that the curves disjoint from $\Gamma_0$ in the open neighborhood also bounds a unique
least area disk in $\Omega$.

\begin{lem} Let $\beta\in V$ with $\beta\cap\Gamma_0= \emptyset$. Then there exist a unique least area disk
spanning $\beta$.
\end{lem}

\begin{pf}
Since $\beta\in V$ is disjoint from $\Gamma_0$, we can find a path $\Gamma:(-\epsilon,\epsilon)\rightarrow
V$, such that $\{\Gamma_t\}$ foliates a neighborhood of $\Gamma_0$ in $\partial \Omega$, and $\beta$ is one
of the leaves, i.e. $\beta=\Gamma_s$ for some $s\in(-\epsilon,\epsilon)$. Then the proofs of the previous
section implies that $\Sigma_\beta=\Sigma_s$ and $\{\Sigma_t\}$ also gives a foliation of a neighborhood of
$\Sigma_\beta$ by least area disks. Then proof of Lemma 5.6 implies that $\Sigma_\beta$ is the unique least
area disk spanning $\beta$.
\end{pf}

Now, if we can show same result for the curves in $V$ intersecting $\Gamma_0$, then we are done.
Unfortunately, we cannot do that, but we will bypass this by going to a smaller neighborhood.

\begin{lem} There exist a neighborhood $V'\subset V$ of $\Gamma_0$ such that for any $\Gamma_0'\in V'$, there exist a
unique least area disk with boundary $\Gamma_0'$.
\end{lem}

\begin{pf}
Let $V'\subset V$ be an open neighborhood containing $\Gamma_0$ such that there exist disjoint two curves
$\beta_1, \beta_2 \in V$ with $\beta_1$ and $\beta_2$ are both disjoint from $\Gamma_0$ and $\Gamma_0'$, for
any $\Gamma_0'\in V'$. We also assume that if $B\subset\partial \Omega$ is the annulus bounded by $\beta_1$
and $\beta_2$, $\Gamma_0, \Gamma_0'$ are contained in $B$. To see the existence of such a neighborhood, one
can fix two curves in $V$ disjoint from $\Gamma_0$, and lying in the opposite sides of $\Gamma_0$ in
$\partial \Omega$. Then suitable complements of these two curves in $V$ will give us the desired neighborhood
of $\Gamma_0$.

Now, fix $\Gamma_0'\in V'$. By the assumption on $V'$, there are two curves $\beta_1,\beta_2$ disjoint from
both $\Gamma_0, \Gamma_0'$ and bounding the annulus $B$ in $\partial \Omega$ such that
$\Gamma_0,\Gamma_0'\subset B\subset \partial \Omega$. Then, we can find two paths
$\Gamma,\Gamma':[-\epsilon,\epsilon]\rightarrow V$ with $\{\Gamma_t\}$, $\{\Gamma_t'\}$ foliates $B$ such
that $\Gamma(\epsilon)={\Gamma'}(\epsilon)=\beta_1$ , $\Gamma({-\epsilon})={\Gamma'}({-\epsilon})=\beta_2$,
and $\Gamma(0)=\Gamma_0$, $\Gamma'(0)=\Gamma_0'$.

By Lemma 5.5, we know that $\{\Gamma_t\}$ induces $\{\Sigma_t\}$ family of embedded least area disks spanning
$\{\Gamma_t\}$. Moreover, these least area disks are unique with the given boundary, and leaves of the
foliation in the neighborhood of $\Sigma_0$.

Now, consider the preimage of the path $\Gamma'$ under the homeomorphism $\pi_1|_U: U \rightarrow V$. This will give us
a path $\Sigma' \subset U\subset M$, which is a continuous family of minimal disks, say $\{\Sigma_t'\}$. We claim that
this is also a family of embedded least area disks inducing a foliated neighborhood of $\Sigma_0'$. By previous
paragraphs, we know that $\Sigma_\epsilon$ and the $\Sigma_{-\epsilon}$ are the unique least area disks with boundary
$\beta_1$ and $\beta_2$, respectively. This means $\Sigma_{\pm\epsilon}'=\Sigma_{\pm\epsilon}$. So, the family
$\{\Sigma_t'\}$ has embedded least area disks $\Sigma_{\pm\epsilon}'$. Then by slight modification of the proof of
Lemma 5.5 imply that $\{\Sigma_t'\}$ is a family of embedded least area disk inducing a foliation of a neighborhood of
of $\Sigma_0'$. By Lemma 5.6, $\Sigma_0'$ is the unique least area disk spanning $\Gamma_0'$.
\end{pf}

So, we got the following theorem.

\begin{thm}
Let $\Omega$ be a $C^3$-smooth mean convex domain in $\BR^3$, and $A= \{\alpha \in C^3(S^1,\partial \Omega) \
| \ \alpha$ embedding$\}$. Then there exist an open dense subset $A'\subset A$ in $C^3$ topology, such that
for any $\Gamma\in A'$, there exist a unique least area disk with boundary $\Gamma$.
\end{thm}

\begin{pf}
The set of regular values of Fredholm map, say $\widehat{A}$, is the whole set except a set of first category
by Sard-Smale theorem. So, the regular curves are dense in $A$. By above lemmas, for any regular curve
$\Gamma_0$, there exist an open neighborhood $V_{\Gamma_0}'\subset A$ which the uniqueness result holds. So,
$A'=\bigcup_{\Gamma \in\hat{A}} V_\Gamma'$ is an open dense subset of $A$ with the desired properties.

\end{pf}

\section{Concluding Remarks}

Different versions of the results in the first part of the paper has been proved by Meeks and Yau in [MY2] by
using differential geometry techniques, by Brian White in [Wh] by using geometric measure theory methods, and
by Fanghua Lin in [Li] by using global analysis methods. Here, we reformulate those results and extend it to
more general class of mean convex domains. Our approach is topological, and seems more natural to the
question. On the other hand, similar results for hyperbolic space has been proved by Michael Anderson in
[An2], by using geometric measure theory techniques.

In this paper, we are trying to promote the idea of {\em mean convex hulls}. These objects are naturally defined for
any simple closed extreme curve, and any asymptotic curve in $\Si$. They are living in convex hulls of the curve, and
have piecewise smooth boundary. As a corollary to the mean convex hulls in $\BH^3$, we assign a {\em mean convex core}
to any quasi-Fuchsian hyperbolic $3$-manifold capturing any minimal surface homotopy equivalent to the manifold.

In the second part, we give a generic uniqueness result, mainly by adapting the techniques of [Co3]. There has been
different types of generic uniqueness results for the curves in $\BR^3$, see [Tr]. This is a new generic uniqueness
result, and says that if you have a mean convex domain with $C^3$-smooth boundary in $\BR^3$, then simple closed curves
in the boundary generically bounds a unique least area disk in $\BR^3$.


\begin{thebibliography}{MSY}


\bibitem[An1]{An1} M. Anderson, {\em Complete minimal varieties in hyperbolic space}, Invent. Math. {\bf 69}, (1982) 477--494.


\bibitem[An2]{An2} M. Anderson, {\em Complete minimal hypersurfaces in hyperbolic n-manifolds}, Comment. Math. Helv. {\bf 58}, (1983)
264--290.

\bibitem[Co1]{Co1} B. Coskunuzer, {\em Uniform 1-cochains and Genuine Laminations}, to appear in Topology.

\bibitem[Co2]{Co2} B. Coskunuzer, {\em Minimal Planes in Hyperbolic Space}, Comm. Anal. Geom. {\bf 12}, (2004) 821--836.

\bibitem[Co3]{Co3} B. Coskunuzer, {\em Generic Uniqueness of Least Area Planes in Hyperbolic Space}, eprint; math.GT/0408066

\bibitem[Ga]{Ga} D. Gabai, {\em On the geometric and topological rigidity of hyperbolic $3$-manifolds}, J. Amer. Math. Soc. {\bf 10},
(1997) 37--74.

\bibitem[HS]{HS} J. Hass and P. Scott, {\em The Existence of Least Area Surfaces in 3-manifolds}, Trans. AMS {\bf 310}, (1988) 87--114.

\bibitem[La]{La} S. Lang, {\em Real analysis}, Addison-Wesley, MA, (1983).

\bibitem[Li]{Li} F.H. Lin, {\em Plateau's problem for $H$-convex curves}, Manuscripta Math. {\bf 58}, (1987) 497--511.

\bibitem[MY1]{MY1} W. Meeks and S.T. Yau, {\em The classical Plateau problem and the topology of three manifolds}, Topology {\bf 21},
(1982) 409--442.

\bibitem[MY2]{MY2} W. Meeks and S.T. Yau, {\em The existence of embedded minimal surfaces and the problem of uniqueness}, Math. Z.
{\bf 179}, (1982) 151--168.

\bibitem[Sm]{Sm} S. Smale, {\em An infinite dimensional version of Sard's Theorem}, Amer. J. Math. {\bf 87}, (1965) 861--866.

\bibitem[TT]{TT} F. Tomi and A.J. Tromba, {\em Extreme curves bound embedded minimal surfaces of the type of the disc}, Math. Z. {\bf
158}, (1978) 137--145.

\bibitem[Tr]{Tr} A. J. Tromba, {\em The set of curves of uniqueness for Plateau's problem has a dense interior}, Geometry
and topology, Lecture Notes in Math., Vol. 597, 696--706, Springer, Berlin, (1977).

\bibitem[Wh]{Wh} B. White, {\em On the topological type of minimal submanifolds}, Topology {\bf 31}, (1992) 445--448.

\end{thebibliography}
\end{document}